\documentclass[12pt]{amsart}
\def\Image{\operatorname{Im}}
\def\Int{\operatorname{Int}}
\def\Min{\operatorname{Min}}
\def\N{{\Bbb N}}

\newtheorem{Theorem}{Theorem}[section]
\newtheorem{Corollary}[Theorem]{Corollary}
\newtheorem{Lemma}[Theorem]{Lemma}

\theoremstyle{definition}

\theoremstyle{remark}

\begin{document}
\sloppy
\title{On non-locally connected boundaries of CAT(0) spaces}
\author{Tetsuya Hosaka} 
\address{Department of Mathematics, Utsunomiya University, 
Utsunomiya, 321-8505, Japan}
\date{April 12, 2006}
\email{hosaka@cc.utsunomiya-u.ac.jp}
\keywords{the boundary of a CAT(0) space}
\subjclass[2000]{57M07}
%\thanks{
%Partly supported by the Grant-in-Aid for Young Scientists (B), 
%The Ministry of Education, Culture, Sports, Science and Technology, Japan.
%(No.~15740029).}
\maketitle
\begin{abstract}
In this paper, 
we study CAT(0) spaces with non-locally connected boundary.
We give some condition of a CAT(0) space 
whose boundary is not locally connected.
\end{abstract}

%%%%%%%%%%%%%
% Section 1 %
%%%%%%%%%%%%%
\section{Introduction and preliminaries}

In this paper, 
we study proper CAT(0) spaces with non-locally connected boundary.
A metric space $X$ is said to be {\it proper} 
if every closed metric ball is compact.
Definitions and basic properties of CAT(0) spaces and their boundaries 
are found in \cite{BH}.

Let $X$ be a proper CAT(0) space and let $\gamma$ be an isometry of $X$.
The {\it translation length} of $\gamma$ is the number 
$|\gamma|:=\inf\{d(x,\gamma x)\,|\,x\in X\}$, and 
the {\it minimal set} of $\gamma$ is defined as 
$\Min(\gamma)=\{x\in X\,|\,d(x,\gamma x)=|\gamma|\}$.
An isometry $\gamma$ of $X$ is said to be {\it hyperbolic}, 
if $\Min(\gamma)\neq\emptyset$ and $|\gamma|>0$ (cf.\ \cite{BH}).
For a hyperbolic isometry $\gamma$ of 
a proper CAT(0) space $X$, 
$\gamma^\infty$ is the limit point of the boundary $\partial X$ 
to which the sequence $\{\gamma^i x_0\}_i$ converges, 
where $x_0$ is a point of $X$.

In this paper, 
we define a {\it reflection} of a geodesic space as follows:
An isometry $r$ of a geodesic space $X$ is called a 
{\it reflection} of $X$, if 
\begin{enumerate}
\item[(1)] $r^2$ is the identity of $X$, 
\item[(2)] $X\setminus F_r$ has exactly 
two convex connected components $X^+_r$ and $X^-_r$ and 
\item[(3)] $rX^+_r=X^-_r$, 
\end{enumerate}
where $F_r$ is the fixed-points set of $r$.
We note that ``reflections'' in this paper need not 
satisfy the condition (4) $\Int F_r=\emptyset$ in \cite{Ho}.

A CAT(0) space $X$ is said to be {\it almost extendible}, 
if there exists a constant $M>0$ such that 
for each pair of points $x,y\in X$, 
there is a geodesic ray $\zeta:[0,\infty)\rightarrow X$ 
such that $\zeta(0)=x$ and $\zeta$ passes within $M$ of $y$. 
In \cite{O}, Ontaneda has proved that 
a CAT(0) space on which some group acts geometrically 
(i.e.\ properly and cocompactly by isometries)
is almost extendible.

In \cite{MR} and \cite{MRT}, 
Mihalik, Ruane and Tschantz have proved some nice results 
about CAT(0) groups with (non-)locally connected boundary.

The purpose of this paper is to prove the following theorem.

\begin{Theorem}\label{Thm}
Let $X$ be a proper and almost extendible CAT(0) space, 
let $\gamma$ be a hyperbolic isometry of $X$ and 
let $r$ be a reflection of $X$.
If 
\begin{enumerate}
\item[(1)] $\gamma^{\infty}\not\in\partial F_r$,
\item[(2)] $\gamma(\partial F_r)\subset\partial F_r$ and
\item[(3)] $\Min(\gamma)\cap F_r=\emptyset$, 
\end{enumerate}
then the boundary $\partial X$ of $X$ is not locally connected.
\end{Theorem}

%%%%%%%%%%%%%
% Section 2 %
%%%%%%%%%%%%%
\section{Topology of the boundary of a CAT(0) space}

In this section, 
we recall topology of the boundary of a CAT(0) space.

Let $X$ be a proper CAT(0) space and $x_0 \in X$.
The {\it boundary of $X$ with respect to $x_0$}, 
denoted by $\partial_{x_0}X$, is defined as 
the set of all geodesic rays issuing from $x_0$. 
Then the topology on $X \cup \partial_{x_0}X$ 
is defined by the following conditions: 
\begin{enumerate}
\item[(1)] $X$ is an open subspace of $X \cup \partial_{x_0}X$. 
\item[(2)] For $\alpha \in \partial_{x_0}X$ and $R, \epsilon >0$, let
$$ U_{x_0}(\alpha;R,\epsilon)=
\{ x \in X \cup \partial_{x_0} X \,|\, 
x \not\in B(x_0,R),\ d(\alpha(R),\xi_{x}(R))<\epsilon \}, $$
where $\xi_{x}:[0,d(x_0,x)]\rightarrow X$ is the geodesic from $x_0$ to $x$
($\xi_{x}=x$ if $x \in \partial_{x_0} X$).
Then for each $\epsilon_0>0$, 
the set 
$$\{U_{x_0}(\alpha;R,\epsilon_0)\,|\, R>0\}$$ 
is a neighborhood basis for $\alpha$ in $X \cup \partial_{x_0}X$. 
\end{enumerate}
This is called the {\it cone topology} on $X \cup \partial_{x_0} X$.
It is known that 
$X \cup \partial_{x_0} X$ is 
a metrizable compactification of $X$ (\cite{BH}, \cite{GH}).

Here the following lemma is known.

\begin{Lemma}\label{Lem}
Let $X$ be a proper CAT(0) space and let $x_0\in X$.
For $\alpha \in \partial_{x_0}X$ and $R, \epsilon >0$, let 
$$ U'_{x_0}(\alpha;R,\epsilon)=
\{ x \in X \cup \partial_{x_0} X \,|\, 
x \not\in B(x_0,R),\ d(\alpha(R),\Image \xi_{x})<\epsilon \}, $$
where $\xi_{x}:[0,d(x_0,x)]\rightarrow X$ is the geodesic from $x_0$ to $x$
($\xi_{x}=x$ if $x \in \partial_{x_0} X$).
Then for each $\epsilon_0>0$, 
the set 
$$\{U'_{x_0}(\alpha;R,\epsilon_0)\,|\, R>0\}$$ 
is also a neighborhood basis for $\alpha$ in $X \cup \partial_{x_0}X$.
\end{Lemma}

Let $X$ be a proper CAT(0) space.
The {\it asymptotic relation} is an equivalence relation 
in the set of all geodesic rays in $X$.
The {\it boundary of} $X$, denoted by $\partial X$, 
is defined as the set of asymptotic equivalence classes of geodesic rays.
The equivalence class of a geodesic ray $\xi$ is denoted by $\xi(\infty)$.
For each $x_0 \in X$ and each $\alpha \in \partial X$, 
there exists a unique element $\xi \in \partial_{x_0}X$ 
with $\xi(\infty)=\alpha$.
Thus we may identify $\partial X$ with $\partial_{x_0}X$ for each $x_0 \in X$
(\cite{BH}, \cite{GH}).

%%%%%%%%%%%%%
% Section 3 %
%%%%%%%%%%%%%
\section{Proof of the theorem}

We prove Theorem~\ref{Thm}.

\begin{proof}[Proof of Theorem~\ref{Thm}]
Let $X$ be a proper and almost extendible CAT(0) space, 
let $\gamma$ be a hyperbolic isometry of $X$ and 
let $r$ be a reflection of $X$ such that 
\begin{enumerate}
\item[(1)] $\gamma^{\infty}\not\in\partial F_r$,
\item[(2)] $\gamma(\partial F_r)\subset\partial F_r$ and
\item[(3)] $\Min(\gamma)\cap F_r=\emptyset$.
\end{enumerate}

Since $X$ is almost extendible, 
there exists a constant $M>0$ such that 
for each pair of points $x,y\in X$, 
there is a geodesic ray $\zeta:[0,\infty)\rightarrow X$ 
such that $\zeta(0)=x$ and $\zeta$ passes within $M$ of $y$. 
By (1), $\gamma^{\infty}\not\in\partial F_r$.
Since $\partial F_r$ is a closed set in $\partial X$, 
there exist $R>0$ and $\epsilon>0$ 
such that $U'_{x_0}(\gamma^{\infty};R,\epsilon)\cap \partial F_r=\emptyset$ 
by Lemma~\ref{Lem}.

Let $x_0\in\Min(\gamma)$ and 
let $\xi:[0,\infty)\rightarrow X$ be the geodesic ray in $X$ 
such that $\xi(0)=x_0$ and $\xi(\infty)=\gamma^{\infty}$.
Then $\Image\xi\subset\Min(\gamma)$.
Since $\Min(\gamma)\cap F_r=\emptyset$ by (3) and 
$\xi(\infty)=\gamma^{\infty}\not\in\partial F_r$ by (1), 
there exists a number $K>0$ such that $d(\xi(K),F_r)>M$.
Let $N=d(x_0,rx_0)$.
For an enough large number $i_0\in\N$, 
$$U'_{x_0}(\gamma^{\infty};i_0|\gamma|,N+K+M)
\subset U'_{x_0}(\gamma^{\infty};R,\epsilon).$$

We prove that 
$U'_{x_0}(\gamma^{\infty};i|\gamma|,N+K+M)\cap\partial X$ 
is not connected for any $i_0\le i\in\N$.
This implies that $\partial X$ is not locally connected, because 
$\{U'_{x_0}(\gamma^{\infty};i|\gamma|,N+K+M)\cap\partial X\,|\,i\in\N, i\ge i_0\}$ 
is a neighborhood basis of $\gamma^{\infty}$ in $\partial X$.

Let $i\in\N$ such that $i\ge i_0$.
Then $\gamma^ir\gamma^{-i}$ is a reflection of $X$ and 
$F_{\gamma^ir\gamma^{-i}}=\gamma^iF_r$.
Here by (3), 
$$ F_{\gamma^ir\gamma^{-i}}\cap\Min(\gamma)=
\gamma^iF_r\cap\gamma^i\Min(\gamma)=
\gamma^i(F_r\cap\Min(\gamma))=\emptyset.$$
Let $X\setminus F_{\gamma^ir\gamma^{-i}}=X_i^+\cup X_i^-$, 
where $X_i^+$ and $X_i^-$ are convex connected components, 
and $x_0\in X_i^+$.
We consider the geodesic ray $\gamma^ir\xi$ such that 
$\gamma^ir\xi(0)=\gamma^irx_0$ and 
$\gamma^ir\xi(\infty)=\gamma^ir\gamma^{\infty}$.
Since $x_0\in X_i^+$, 
$\Image\xi\subset X_i^+$ and $\Image\gamma^ir\xi\subset X_i^-$.
Hence $\gamma^ir\xi(K)\in X_i^-$.
Here 
\begin{align*}
&d(\gamma^irx_0,\gamma^ir\xi(K))=d(x_0,\xi(K))=K \ \text{and} \\
&d(\gamma^ir\xi(K),\gamma^irF_r)=d(\xi(K),F_r)>M.
\end{align*}
By the definition of the number $M$, 
there exists a geodesic ray $\zeta_i:[0,\infty)\rightarrow X$ 
such that $\zeta_i(0)=x_0$ and 
$\zeta_i$ passes within $M$ of $\gamma^ir\xi(K)$.
Since $d(\gamma^ir\xi(K),\gamma^irF_r)>M$, 
$\zeta_i(\infty)\in \partial X_i^-$.
Because if $\zeta_i(\infty)\in \partial X\setminus\partial X_i^-=
\partial (X_i^+\cup\gamma^iF_r)$ 
then $\Image\zeta_i\subset X_i^+\cup \gamma^iF_r$, 
since $X_i^+\cup \gamma^iF_r$ is convex and $\zeta_i(0)=x_0\in X_i^+$.
Then 
\begin{align*}
d(\gamma^ix_0,\Image\zeta_i)&\le 
d(\gamma^ix_0,\gamma^irx_0)+d(\gamma^irx_0,\gamma^ir\xi(K))
+d(\gamma^ir\xi(K),\Image\zeta_i) \\
&\le d(x_0,rx_0)+d(x_0,\xi(K))+M \\
&=N+K+M.
\end{align*}
We note that $\xi(\infty)=\gamma^\infty$ and $\gamma^ix_0=\xi(i|\gamma|)$, 
since $x_0\in\Min(\gamma)$.
Hence $\zeta_i\in U'_{x_0}(\gamma^\infty;i|\gamma|,N+K+M)$.

Now we show that there does not exist 
a path from $\gamma^\infty$ to $\zeta_i(\infty)$ 
in $U'_{x_0}(\gamma^\infty;i|\gamma|,N+K+M)\cap\partial X$.
Since $\gamma^\infty\in\partial X_i^+$ and 
$\zeta_i(\infty)\in\partial X_i^-$, 
such pass must intersect with $\partial F_{\gamma^ir\gamma^{-i}}$.
Here 
$$\partial F_{\gamma^ir\gamma^{-i}}=\partial(\gamma^iF_r)
=\gamma^i(\partial F_r)\subset\partial F_r,$$
by (2).
We note that 
$U'_{x_0}(\gamma^{\infty};R,\epsilon)\cap \partial F_r=\emptyset$ and 
\begin{align*}
U'_{x_0}(\gamma^{\infty};i|\gamma|,N+K+M)&\subset 
U'_{x_0}(\gamma^{\infty};i_0|\gamma|,N+K+M) \\
&\subset U'_{x_0}(\gamma^{\infty};R,\epsilon).
\end{align*}
Hence 
$$U'_{x_0}(\gamma^{\infty};i|\gamma|,N+K+M)\cap \partial F_{\gamma^ir\gamma^{-i}}
=\emptyset.$$
Thus there does not exist a path between $\gamma^\infty$ and $\zeta_i(\infty)$ 
in $U'_{x_0}(\gamma^\infty;i|\gamma|,N+K+M)\cap\partial X$.

Therefore $\partial X$ is not locally connected.
\end{proof}

%%%%%%%%%%%%%
% Section 3 %
%%%%%%%%%%%%%
\section{Remark}

Every CAT(0) space on which 
some group acts geometrically (i.e.\ properly and cocompactly by isometries) 
is proper (\cite[p.132]{BH}) and almost extendible (\cite{O}).

In \cite{Ru}, Ruane has proved that $\partial \Min(\gamma)$ is 
the fixed-points set of $\gamma$ in $\partial X$, i.e., 
$$\partial \Min(\gamma)=\{\alpha\in\partial X\,|\, \gamma\alpha=\alpha \}.$$
Hence, for example, 
if $\partial F_r \subset \partial\Min(\gamma)$ 
then $\gamma(\partial F_r)=\partial F_r$ and 
the condition (2) in Theorem~\ref{Thm} holds.

A Coxeter system $(W,S)$ defines a {\it Davis complex} $\Sigma(W,S)$ 
which is a CAT(0) space (\cite{D1} and \cite{M}).
Then the Coxeter group $W$ acts geometrically on $\Sigma(W,S)$ and 
each $s\in S$ is a reflection of $\Sigma(W,S)$.

For example, as an application of Theorem~\ref{Thm}, 
we can obtain the following corollary.

\begin{Corollary}\label{Cor}
Let $(W,S)$ be a right-angled Coxeter system and 
let $\Sigma(W,S)$ be the Davis complex of $(W,S)$.
Suppose that 
there exist $s_0,s_1,u_0\in S$ such that 
\begin{enumerate}
\item[(1)] $o(s_0s_1)=\infty$, 
\item[(2)] $o(s_0u_0)=\infty$ and 
\item[(3)] $s_0t=ts_0$ and $s_1t=ts_1$ for each $t\in \tilde{T}$,
\end{enumerate}
where $T=\{t\in S\,|\,tu_0=u_0t\}$ and $\tilde{T}$ is the subset of $S$ 
such that $W_{\tilde{T}}$ is 
the minimum parabolic subgroup of finite index in $W_T$.
Then the boundary $\partial\Sigma(W,S)$ is not locally connected.
\end{Corollary}

\begin{proof}
Let $\gamma=s_0s_1$ and $r=u_0$.
Then $\gamma$ is a hyperbolic isometry of $\Sigma(W,S)$ by (1), 
$r$ is a reflection of $\Sigma(W,S)$ and 
$\partial F_r=\partial\Sigma(W_{\tilde{T}},\tilde{T})$.
Here by (3), 
$$\gamma(\partial F_r)=(s_0s_1)\partial\Sigma(W_{\tilde{T}},\tilde{T})
=\partial\Sigma(W_{\tilde{T}},\tilde{T})=\partial F_r.$$
Also 
$$\gamma^\infty=(s_0s_1)^\infty\not\in
\partial\Sigma(W_{\tilde{T}},\tilde{T})=\partial F_r,$$
and $\Min(\gamma)\cap F_r=\emptyset$ by (2).
Thus the conditions in Theorem~\ref{Thm} hold, 
and $\partial\Sigma(W,S)$ is not locally connected.
\end{proof}

Corollary~\ref{Cor} is a special case of Theorem~3.2 in \cite{MRT}.
We can also obtain Corollary~\ref{Cor} from Theorem~3.2 in \cite{MRT}.

%%%%%%%%%%%%%%%%%%%%%%%%%%%%%%%%%%%%%
%             REFERENCES            %
%%%%%%%%%%%%%%%%%%%%%%%%%%%%%%%%%%%%%
%

%
\end{document}